\documentclass[12pt, reqno]{article}
\usepackage{helvet}
\usepackage{amsfonts}
\usepackage{latexsym}
\usepackage{amsmath}
\usepackage{amsthm,amssymb,mathrsfs,amstext}
\font\msbm=msbm10

\numberwithin{equation}{section}
 
 \def\cS{\mathcal{S}}

 \def\O{\mathcal{O}}

\theoremstyle{plain}
\newtheorem{Theorem}{Theorem}[section]
\newtheorem{lemma}[Theorem]{Lemma}
\newtheorem{corollary}[Theorem]{Corollary}
\newtheorem{proposition}[Theorem]{Proposition}

\newtheorem{remark}{Remark}[section]

\newtheorem{definition}{Definition}[section]

\title{Coorbit Space Theory for Quasi-Banach Spaces}
\author{Holger Rauhut\\
NuHAG, Faculty of Mathematics, University of Vienna\\
Nordbergstrasse 15, A-1090 Wien, Austria\\
rauhut@ma.tum.de}
\date{}

\def\mathbb#1{\hbox{\msbm{#1}}}
\newcommand{\N}{{\mathbb{N}}}
\newcommand{\R}{{\mathbb{R}}}
\newcommand{\Z}{{\mathbb{Z}}}

\newcommand{\TT}{{\mathbb{T}}}
\newcommand{\BB}{{\mathbb{B}}}
\newcommand{\CC}{{\cal{C}}}
\newcommand{\A}{{\mathbb{A}}}
\newcommand{\HH}{{\mathbb{H}}}

\newcommand{\bS}{{\mathbb{S}}}

\newcommand{\G}{{\cal{G}}}

\renewcommand{\H}{{\cal{H}}}
\newcommand{\D}{{\cal{D}}}
\newcommand{\DD}{{\mathbb{D}}}
\newcommand{\Co}{{\mathsf{Co}}}

\newcommand{\on}{{|\!|\!|}}

\newcommand{\STFT}{\operatorname{STFT}}

\newcommand{\beq}{\begin{eqnarray}}
\newcommand{\eeq}{\end{eqnarray}}
\newcommand{\beqn}{\begin{eqnarray*}}
\newcommand{\eeqn}{\end{eqnarray*}}

\newcommand{\ol}{\overline}
\newcommand{\supp}{\operatorname{supp}}

\renewcommand{\qed}{\rule{2.5mm}{2.5mm}}
\newenvironment{Proof}{\noindent
{\bf\underline{Proof:} }}
{\hspace*{\fill}\qed\vskip1em}

\begin{document}
\maketitle
\begin{abstract}
We generalize the classical coorbit space theory developed by Feich\-tinger
and Gr\"ochenig to quasi-Banach spaces. 
As a main result we provide atomic 
decompositions for coorbit spaces defined with respect to quasi-Banach spaces. 
These atomic decompositions are used to prove
fast convergence rates of best $n$-term approximation schemes. 
We apply the abstract
theory to time-frequency analysis of modulation spaces 
$M^{p,q}_m$, $0< p,q \leq \infty$.
\end{abstract}




\section{Introduction}

Coorbit space theory was originally developed by Feichtinger and Gr\"ochenig
\cite{FG1,FG2,FG3,gro} in the late 1980's with the aim 
to provide a unified and group-theoretical approach to function spaces and
their atomic decompositions. In particular, this theory covers the homogeneous
Besov and Triebel-Lizorkin spaces and their wavelet-type atomic decompositions,
as well as the modulation spaces and their Gabor-type decompositions.
Recently, there has been some activity to provide 
generalizations to other settings than the classical one of 
integrable group representations \cite{DST,DST1,FR,hr,hr2}.

All the approaches done so far cover only the case of Banach spaces.
For certain applications such as non-linear approximation, however, 
it is useful to consider also the case of quasi-Banach spaces. For instance 
this would allow to describe also modulation spaces
$M^{p,q}_m$ with $p<1$ or $q<1$, or Hardy spaces $H^p$ with $p<1$, as 
coorbit spaces. In \cite{FG1} it is remarked 
that such an extension of
coorbit space theory to quasi-Banach spaces would be interesting, 
but it seems that 
nothing concrete has been done since then.

So this paper deals with
such an extension of the classical coorbit space theory. 
Our starting point is an integrable representation $\pi$ 
of some locally compact
group $\G$ on some Hilbert space $\H$. Associated 
to $\pi$ is the abstract wavelet transform 
$V_g f(x) = \langle f,\pi(x)g \rangle$.
The crucial ingredient in coorbit space theory is the reproducing formula 
for $V_g$, see (\ref{rep_form_H}), which uses the group 
convolution on $\G$. Thus, it is essential to have a 
convolution relations for certain quasi-Banach spaces $Y$ on $\G$.
Unfortunately, even for the natural choice $Y=L^p(\G)$, $0<p<1$, no 
convolution relation is available. In order to overcome
this problem we work with Wiener amalgam spaces $W(L^\infty,Y)$
with local component $L^\infty$ instead of $Y$ itself. 
Convolution relations for such spaces, where $Y$ is allowed 
to be a quasi-Banach space, were shown recently by the author in
\cite{Rau_Wiener}.

Under some technical assumption on the representation, 
the coorbit spaces $\CC(Y)$ are defined as retract of the Wiener amalgam space
$W(L^\infty,Y)$ via the abstract wavelet transform, i.e.,
$\CC(Y) \,=\, \{f, V_g f \in W(L^\infty,Y)\}$. 
We will prove that $\CC(Y)$ is indeed a quasi-Banach 
space that is independent of the choice of $g$. Moreover,
analogously as in \cite{FG1,FG2,gro} we will provide atomic
decompositions of $\CC(Y)$ of the form $\{\pi(x_i) g\}_{i \in I}$, where
$(x_i)_{i\in I}$ is a suitable point set in the group.

Our results are applicable to time-frequency analysis on modulation spaces 
$M^{p,q}_m$, $0<p,q\leq \infty$, introduced by Feichtinger \cite{Fei_mod}, see also
\cite{Triebel,GS} for the case $p,q<1$. 
Hereby, we improve or give alternative proofs 
to some of the results
of Galperin and Samarah in \cite{GS}. In particular, we show that regular Gabor 
frame expansions with a Schwartz class window automatically extend from $L^2(\R^d)$
to all modulation spaces $M^{p,q}_m$, $0<p,q\leq\infty$ (with moderate 
weight functions $m$ of polynomial growth).

We remark that the abstract theory applies also to homogeneous (weighted)
Besov spaces $\dot{B}_{p,q}^s$ and Triebel-Lizorkin spaces $\dot{F}_{p,q}^s$,
$0<p,q\leq \infty$. We postpone a detailed discussion to a subsequent contribution.

The paper is organized as follows. In Section 2 we introduce some notation
and prerequisites.
Section 3 recalls recent results from \cite{Rau_Wiener} 
on Wiener amalgam spaces with respect
to quasi-Banach spaces including convolution relations. 
Then in Section 4 we introduce the coorbit spaces and show their
basic properties. The atomic decompositions of the coorbit spaces 
will be provided
in Section 5 and Section 6 deals with the question whether the coorbit space
admits also a characterization by $Y$ itself rather than by $W(L^\infty,Y)$.
Section 7 investigates the approximation rates of the best $n$-term approximation
with elements of the atomic decomposition. Here, Lorentz spaces play a key
role. Finally, we apply our abstract results to time-frequency analysis
of modulation spaces in Section 8.

\section{Prerequisites}

Let $\G$ be a locally compact group with identity $e$. 
Integration on $\G$ will always
be with respect to the left Haar measure. We denote by $L_x F(y) = F(x^{-1}y)$ and
$R_x F(y) = F(yx)$, $x,y \in \G$, the left and right translation operators.
Furthermore, let $\Delta$ be the Haar-module on $\G$. 
For a Radon measure $\mu$ we introduce the operator 
$(A_x \mu)(k) \:= \mu(R_x k)$, $x \in \G$,
for a continuous function $k$ with compact support. 
We may identify 
a function $F \in L^1$
with a measure $\mu_F \in M$ by $\mu_F(k) = \int F(x) k(x) dx$. 
Then it clearly holds $A_x F = \Delta(x^{-1}) R_{x^{-1}} F$.
Further, we define the
involutions $F^\vee(x) = F(x^{-1})$, 
$F^\nabla(x) = \overline{F(x^{-1})}$, $
F^*(x) = \Delta(x^{-1}) \overline{F(x^{-1})}$. 

A quasi-norm $\|\cdot\|$ 
on some linear space $Y$ is defined in the same way as a norm, with
the only difference that the triangle inequality
is replaced by $\|f+g\| \leq C(\|f\|+\|g\|)$ with some constant 
$C\geq 1$. 
It is well-known, 
see e.g. \cite[p.~20]{DL} or \cite{Kalton}, that 
there exists an equivalent quasi-norm $\|\cdot|Y\|$  
on $Y$ and an exponent $p$ with $0 < p \leq 1$ such that $\|\cdot|Y\|$
satisfies the $p$-triangle inequality, i.e., 
$\|f+g|Y\|^p \leq \|f|Y\|^p + \|g|Y\|^p$. 
We can choose $p=1$ 
if and only if $Y$ is a normed space. We always
assume in the sequel that such a $p$-norm 
on $Y$ is chosen and denote it by $\|\cdot|Y\|$. If $Y$ is complete
with respect to the topology defined by the metric $d(f,g) = \|f-g|Y\|^p$
then it is called a quasi-Banach space.

A quasi-Banach space of measurable
functions on $\G$ is called solid
if $F \in Y$, $G$ measurable and satisfying $|G(x)| \leq |F(x)|$ a.e. implies
$G \in Y$ and $\|G|Y\| \leq \|F|Y\|$. 
The Lebesgue spaces $L^p(\G)$, $0<p\leq \infty$, provide natural 
examples of solid quasi-normed spaces on $\G$, 
and the usual quasi-norm in $L^p(\G)$ is a 
$p$-norm if $0<p\leq 1$.
If $w$ is some positive measurable weight function on $\G$ then we further define 
$L^p_w = \{F \mbox{ measurable }, Fw \in L^p\}$ with  
$\|F|L^p_w\| := \|Fw|L^p\|$. A continuous weight $w$ is called submultiplicative
if $w(xy) \leq w(x) w(y)$ for all $x,y \in \G$. Further, another weight $m$ is called
$w$-moderate if $m(xyz) \leq w(x) m(y) w(z)$, $x,y,z \in \G$.
It is easy to see that $L^p_m$
is invariant under left and right translations if $m$ is $w$-moderate.

For a quasi-Banach space $(B,\|\cdot|B\|)$ we denote the quasi-norm of a 
bounded operator $T:B\to B$ by $\on T|B\on$. 
The symbol $A\asymp B$ indicates throughout
the paper that there are constants $C_1, C_2 > 0$ such that $C_1 A \leq B \leq C_2 A$
(independently on other expressions on which $A,B$ might depend). The symbol
$C$ will always denote a generic constant whose precise value might
differ at different occurences.

\section{Wiener Amalgam Spaces}

Let $B$ be one of the spaces $L^\infty(\G), L^1(\G)$ or $M(\G)$, 
the space of complex Radon measures.  
Choose some relatively compact neighborhood $Q$ of $e \in \G$.
We define the control function by
\begin{equation}\label{def_control}
K(F,Q,B)(x) \,:=\, \|(L_x \chi_Q) F|B\|, \quad x\in \G,
\end{equation}
where $F$ is locally contained in $B$, in symbols $F\in B_{loc}$.
Further, let $Y$ be some solid quasi-Banach space of functions
on $\G$ containing the characteristic function of any compact subset of $\G$. 
The {\bf Wiener amalgam space} $W(B,Y)$ is then defined by
\[\label{def_Wiener_space}
W(B,Y) \,:=\, W(B,Y,Q) \,:=\, \{F \in B_{loc},\, K(F,Q,B) \in Y\}
\]
with quasi-norm
\begin{equation}\label{qnormW}
\|F|W(B,Y,Q)\|\,:=\, \|K(F,Q,B)|Y\|.
\end{equation}
This is indeed a $p$-norm with $p$ being the exponent of the quasi-norm of $Y$.
By $W(C_0,Y)$ we denote the closed subspace of $W(L^\infty,Y)$ consisting of
continuous functions.

We also need certain discrete sets in $\G$. 

\begin{definition} Let $X=(x_i)_{i\in I}$ be some discrete set of points in $\G$ and 
$V$, $W$ relatively compact neighborhoods of $e$ in $\G$.
\begin{itemize}\itemsep=-1pt
\item[(a)] $X$ is called $V$-dense if $\G = \bigcup_{i \in I} x_i V$.
\item[(b)] $X$ is called relatively 
separated if for all compact sets
$K \subset \G$ there exists a constant $C_K$ such that
$\sup_{j \in I} \#\{ i\in I,\, x_iK \cap x_jK \neq \emptyset \} \leq C_K. 
$
\item[(c)] $X$ is called $V$-well-spread (or simply well-spread) if it is 
both relatively separated and $V$-dense for
some $V$.
\end{itemize}
\end{definition}
The existence of $V$-well-spread sets for arbitrarily small $V$ is proven in
\cite{Fei_Homog}, see also \cite{hr,hr2} for a generalization.
Given a well-spread family $X=(x_i)_{i\in I}$, a 
relatively compact neighborhood $Q$ of $e \in \G$ 
and  $Y$, we define the sequence space 
\begin{align}\label{def_Ydiscrete}
Y_d \,:=\, Y_d(X) \,:=\, Y_d(X,Q) \,:=&\, \{ (\lambda_i)_{i \in I}, \sum_{i\in I} 
|\lambda_i|\chi_{x_i Q} \in Y\},
\end{align}
with natural norm
$\|(\lambda_i)_{i \in I} | Y_d\| := \|\sum_{i \in I} |\lambda_i| \chi_{x_i Q}|Y\|.$
Hereby, $\chi_{x_i Q}$ denotes the characteristic function of the set $x_i Q$. 
If the quasi-norm of $Y$ is a $p$-norm,
$0 < p \leq 1$, then also $Y_d$ has a $p$-norm.  

The following concept will also be very useful.

\begin{definition}\label{def_IBUPU} Suppose $U$ is a relatively compact
neighborhood of $e \in \G$. A collection of functions 
$\Psi = (\psi_i)_{i\in I}, \psi_i \in C_0(\G)$, is called
bounded uniform partition of unity of size $U$ 
(for short $U$-BUPU) if the 
following conditions are satisfied:
\begin{itemize}\itemsep-1pt
\item[(1)] $0 \leq \psi_i(x) \leq 1$ for all $i\in I$, $x\in \G$,
\item[(2)] $\sum_{i \in I} \psi_i(x) \equiv 1$,
\item[(3)] there 
exists a  well-spread family $(x_i)_{i\in I}$ such that 
$
\supp \psi_i \,\subset\, x_i U.
$     
\end{itemize}
\end{definition}

The construction of BUPU's with respect to arbitrary well-spread sets 
is standard.

We call $W(B,Y,Q)$ right translation invariant 
if the right translations $R_x$ (resp. $A_x$ if $B=M$) 
are bounded operators on $W(L^\infty,Y,Q)$.
In \cite{Rau_Wiener} the following results were shown.

\begin{Theorem}\label{thm_basic} The following statements are equivalent:
\begin{itemize}\itemsep-1pt
\item[(i)] $W(L^\infty,Y) = W(L^\infty,Y,Q)$ is independent of the choice
of the neighborhood $Q$ of $e$ (with equivalent norms for different choices). 
\item[(ii)] $Y_d = Y_d(X,U)$ is independent of the choice of the neighborhood
 $U$ of $e$ (with equivalent norms for different choices).
\item[(iii)] $W(L^\infty,Y) = W(L^\infty,Y,Q)$ is right translation invariant.
\end{itemize} 
If one (and hence all) of these conditions are satisfied then:
\begin{itemize}\itemsep-1pt
\item[(a)] $W(B,Y)=W(B,Y,Q)$ is 
independent of the choice of $Q$. 
\item[(b)] $W(B,Y)$ is right translation invariant.
\item[(c)] $Y_d$ and $W(B,Y)$ are complete.
\item[(d)] The expression
\begin{equation}\label{WYd_norm}
\|F|W(B, Y_d)\| \,:=\, \| (\|F \chi_{x_i Q}|B\|)_{i\in I} | Y_d(X)\|,
\end{equation} 
defines an equivalent quasi-norm on $W(B,Y)$.
\end{itemize}
\end{Theorem}

 
Also the left translation invariance is a useful property. For instance,
it ensures inclusions into weighted $L^\infty$ spaces, see \cite{Rau_Wiener}.

\begin{lemma}\label{lem_Linfty_embed} Let $W(L^\infty,Y)$ be
left translation invariant. Let $r(x):=  \on L_{x^{-1}}|W(L^\infty,Y)\on$ 
and $\tilde{r}(i):= r(x_i)$. Then
\begin{itemize}
\item[(a)] $Y_d(X)$ is continuously embedded into $\ell^\infty_{1/\tilde{r}}$;
\item[(b)] $W(L^\infty,Y)$ is continuously embedded into $L^\infty_{1/r}$.
\end{itemize}
\end{lemma} 

The following criterions for left and right translation invariance
were provided in \cite{Rau_Wiener}.

\begin{lemma}\label{lem_left_trans} 
If $Y$ is left translation invariant then $W(B,Y)$
is left translation invariant and 
$\on L_y |W(B,Y)\on \leq \on L_y | Y \on$. 
\end{lemma}

Recall that $\G$ is called an IN group if there exists a compact neighborhood of $e$
such that $xQ = Qx$ for all $x \in \G$. 

\begin{lemma}\label{lem_IN_trans} Let $Y$ be right translation invariant.
Then also $W(B,Y)$ is right translation invariant. Moreover,
if $\G$ is an IN group then 
\[
\on R_y|W(L^\infty,Y)\on \leq \on R_y|Y\on
\quad\mbox{ and }\quad 
\on A_y|W(M,Y)\on \leq \on R_y|Y\on.
\]
\end{lemma}

We remark that $Y$ does not necessarily need to be translation invariant in 
order $W(L^\infty,Y)$ to be
translation invariant, see \cite{Rau_Wiener} for an example.
The following result for the involution ${ }^\vee$ will also be useful 
later on.

\begin{lemma}\label{lem_IN_inv}
If $\G$ is an IN group then $W(L^\infty,Y^\vee)^\vee = W(L^\infty,Y)$. 
\end{lemma}

The main ingredient for the coorbit space theory with respect to quasi-Banach 
spaces will be
the following convolution relations for Wiener amalgam spaces.

\begin{Theorem}\label{thm_conv} Let $0 < p \leq 1$ be such that the quasi-norm of 
$Y$ satisfies the $p$-triangle inequality and assume that $W(L^\infty,Y)$ is right translation 
invariant.
\begin{itemize}
\item[(a)] Set $w(x):= \on A_x|W(M,Y) \on$. Then we have 
\[
W(M,Y) * W(L^\infty,L^p_w) \hookrightarrow W(L^\infty,Y)
\]
with corresponding estimate for the quasi-norms.
\item[(b)] Set $v(x):= \Delta(x^{-1}) \on R_{x^{-1}} | W(L^\infty,Y)\|$. Then we have
\[
W(L^\infty,Y) * W(L^\infty,L^p_v) \hookrightarrow W(L^\infty,Y)
\]
with corresponding estimate for the quasi-norms. 
\end{itemize} 
\end{Theorem}

%
%



\begin{Theorem}\label{thm_conv_Lp} 
Let $w$ be a submultiplicative weight and $0<p\leq 1$. Then it holds
\[
W(L^\infty,L^p_w) * W(L^\infty,L^p_{w^*})^\vee \hookrightarrow W(L^\infty,L^p_w).
\] 
In particular, if $\G$ is an IN-group then 
$W(L^\infty,L^p_w) * W(L^\infty,L^p_w) \hookrightarrow W(L^\infty,L^p_w)$ with corresponding
quasi-norm estimate.
\end{Theorem}



Further, we will need the following maximal function.
For some relatively compact neighborhood $U$ of $e \in \G$ and a function
$G$ on $\G$ we define the $U$-oscillation by
\[
G_U^\#(x) \,:=\, \sup_{u \in U} |G(ux) - G(x)|.
\]

The following lemma on the $U$-oscillation will be an 
essential tool for deriving the atomic decomposition
for the coorbit spaces defined later on.

\begin{lemma}\label{lem_osc}
\begin{itemize}
\item[(a)] Let $W(L^\infty,Y)$ be left and right translation invariant. If 
$G \in W(C_0,W(L^\infty,Y)^\vee)^\vee \cap W(C_0,Y)$ then
$G_U^\# \in W(C_0,Y)$.  
\item[(b)] 
Let $w$ be a submultiplicative weight function and $0<p< \infty$.
Then $G \in W(C_0,W(L^\infty,L^p_w)^\vee)^\vee\cap W(C_0,L^p_w)$ implies
\[
\lim_{U \to \{e\}} \|G^\#_U|W(L^\infty,L^p_w)\| = 0.
\]
\item[(c)] If $y \in xU$ then $|L_y G - L_xG| \leq L_y G^\#_U$ holds pointwise.
\end{itemize}
\end{lemma}
\begin{Proof}
(a) The left and right translation invariance of $W(L^\infty,Y)$ 
implies by Theorem \ref{thm_basic} that
$W(C_0,W(L^\infty,Y)^\vee,Q)$ and $W(C_0,Y,Q)$
are independent of the choice of the compact neighborhood $Q$ of $e \in \G$.
The control function of $G^\#_U$ can be estimated as follows,
\begin{align}
&K(G^\#_U,Q,L^\infty)(x) \,=\, \sup_{z \in xQ} G^\#_U(z)
\,=\, \sup_{z \in xQ} \sup_{u\in U} |G(uz) -G(z)|\notag\\
&\leq\, \sup_{z \in xQ} \sup_{u\in U} |G(uz)| + \sup_{z \in xQ} |G(z)|
\,=\, \sup_{q \in Q} \sup_{u\in U} |G(uxq)| + K(G,Q,C_0)(z).\notag
\end{align}
Clearly, we have $K(G,Q,C_0) \in Y$ by assumption on $G$.
We further compute the function 
$H(x):= \sup_{q \in Q} \sup_{u\in U} |G(uxq)|$
\begin{align}
H(x) \,&=\, \sup_{q\in Q} \|\chi_U (R_{xq} G)\|_\infty
\,=\, \sup_{q\in Q} \|(R_{(xq)^{-1}} \chi_U)^\vee G^\vee\|_\infty\notag\\
&=\, \sup_{q\in Q} \|L_{(xq)^{-1}} \chi_{U^{-1}} G^\vee\|_\infty
\,=\, \sup_{q\in Q} K(G^\vee,U^{-1},L^\infty)^\vee(xq)\notag\\
&=\, \|\chi_{xQ} K(G^\vee,U^{-1},L^\infty)^\vee\|_\infty
\,=\, K(K(G^\vee,U^{-1},L^\infty)^\vee,Q,L^\infty)(x).\notag
\end{align}
Thus,
\begin{align}
\|H|Y\| \,&=\, \|K(K(G^\vee,U^{-1},L^\infty)^\vee,Q,L^\infty)|Y\|\notag\\
&\leq\,C \|G^\vee|W(L^\infty,W(L^\infty,Y)^\vee)\|.\notag
\end{align}
This implies $G_U^\# \in W(C_0,Y)$.

(b) By part (a) $G_U^\#$ is contained in $W(C_0,L^p_w)$. 
Let $\epsilon > 0$. Since 
$U \subset U_0$ implies $G_U^\# \leq G_{U_0}^\#$ we can find 
a compact set $V \subset \G$
such that 
\[
\int_{\G\setminus V} K(G_U^\#,Q,L^\infty)(x)^p w(x)^p dx \leq \frac{\epsilon}{2}
\]
for all $U \subset U_0$. Since $G$ is uniformly continuous on the compact set 
$VQ$ we can find a neighborhood $U_1 \subset U_0$ of $e$ such that
\[
G_{U_1}^\#(x) \,\leq\, M\,:=\, \frac{\epsilon^{1/p}}{(2|V|)^{1/p}\, \nu} \quad \mbox{ for all } x\in VQ
\]
with $\nu:= \max_{x \in V} w(x)$. This implies
\[
K(G_{U_1}^\#,Q,L^\infty)(x) \,=\, \sup_{z \in xQ} |G_{U_1}^\#(z)| \,\leq\, M \quad \mbox{ for all }x \in V.
\]
Thus, we obtain
\[
\int_V K(G_{U_1}^\#,Q,L^\infty)(x)^p w(x)^p dx
\,\leq\, M^p |V| \nu^p \,=\, \frac{\epsilon}{2}.
\]
Altogether this yields $\|G_{U_1}^\#|W(C_0,L^p_w)\|^p \leq \epsilon$.

(c) This is straightforward (see also Lemma 4.6(iii) in \cite{gro}).
\end{Proof}

\begin{remark}\label{rem_GIN} Let $\G$ be an IN-group. Then it follows from
Lemma \ref{lem_IN_inv} that
\[
W(C_0,W(L^\infty,Y)^\vee)^\vee \,=\, W(C_0,W(L^\infty,Y)) \,=\, W(C_0,Y).
\]
The second equality follows from $K(K(F,Q,L^\infty))\leq K(F,Q^2,L^\infty)$
and the independence of $W(L^\infty,Y,Q)$ of $Q$.
Thus, it suffices to assume $G \in W(C_0,L^p_w)$ in (b) in this case. For 
general groups, however, such a simplification does not seem possible.
\end{remark}

\section{Coorbit Spaces}

Let $\pi$ be an irreducible unitary representation of $\G$
on some Hilbert space $\H$. Then the abstract wavelet transform (voice
transform) is defined as
\[
V_g f(x) \,:=\, \langle f, \pi(x) g\rangle, \qquad f,g \in \H, x\in \G. 
\]
The representation $\pi$ is called square-integrable if there exists a non-zero
$g \in \H$ (called admissible) such that $V_g g \in L^2(\G)$. Then by a 
theorem of Duflo and Moore \cite{DM}
it holds 
\begin{equation}\label{isometry}
\|V_g f|L^2\| \,=\, c_g \|f|\H\| \qquad \mbox{ for all }f \in \H
\end{equation}
with some constant
$c_g$. It can be shown \cite{DM} that $c_g = \|Kg|\H\|$ for some 
uniquely defined
self-adjoint, positive and densely defined operator $K$ (possibly unbounded) whose domain
$\D(K)$ consists of the admissible vectors. Moreover, if $\G$ is unimodular then
$K$ is a multiple of the identity. 

As a consequence of (\ref{isometry}), if $g$ is normalized, i.e., $c_g = \|Kg\|=1$, 
we have the reproducing formula
\begin{equation}\label{rep_form_H}
V_g f \,=\, V_g f * V_g g.
\end{equation}
In order to introduce the coorbit spaces we first need to extend
the definition of the voice transform to a larger space, the ``reservoir''.
To this end let $v$ be some submultiplicative
weight function satisfying $v\geq 1$. We define the following class
of analyzing vectors
$$
\A_v \,:=\, \{g \in \H, V_g g \in L^1_v\}.
$$
We assume that $\A_v$ is non-trivial, i.e., $\pi$ is integrable. This implies
that $\pi$ is also square-integrable. With some fixed $g \in \A_v \setminus \{0\}$
we define
\[
\H_v^1 \,:=\, \{f \in \H, V_g f \in L^1_v\}
\]
with norm $\|f|\H_v^1\|:= \|V_g f|L^1_v\|$.
It can be shown \cite{FG1,hr} that $\H_v^1$ is a $\pi$-invariant Banach space whose
definition does not depend on the choice of $g$. We denote by $(\H_v^1)^\urcorner$ the anti-dual, i.e.,
the space of all bounded conjugate-linear functionals on $\H^1_v$.
An equivalent norm
on $(\H^1_v)^\urcorner$ is given by $\|V_g f|L^\infty_{1/v}\|$.
Denoting by $\langle \cdot,\cdot\rangle$ also the dual pairing
of $(\H^1_v,(\H^1_v)^\urcorner)$ the voice transform 
extends to $(\H^1_v)^\urcorner$,
\[
V_g f(x) \,=\, \langle f,\pi(x) g\rangle, \quad f \in (\H^1_v)^\urcorner, g\in \A_v.
\]
Important properties of the voice transform extend to $(\H^1_v)^\urcorner$ as stated
in the following lemma, see \cite{FG1,FG2,hr}.
\begin{lemma}\label{lem_antidual} Let $g \in \A_v$ with $\|Kg|\H\|=1$.
\begin{itemize} 
\item[(a)] It holds $V_g (\pi(x)f) = L_x V_g f$ 
for all $x \in \G$, $f \in (\H^1_v)^\urcorner$.
\item[(b)] A bounded net $(f_\alpha) \subset (\H^1_v)^\urcorner$ is weak-$*$ convergent
to an element $f \in (\H^1_v)^\urcorner$ if and only if $(V_g f_\alpha)$ converges
to $V_g f$ pointwise.
\item[(c)] The reproducing formula extends to $(\H^1_v)^\urcorner$, i.e.,
\[
V_g f \,=\, V_gf * V_g g \qquad \mbox{for all }f\in (\H^1_v)^\urcorner.
\]
\item[(d)] Conversely, if $F \in L^\infty_{1/w}$ satisfies the reproducing formula
$F=F*V_g g$ then there exists a unique element $f\in (\H^1_v)^\urcorner$ such that
$F= V_g f$.
\end{itemize}
\end{lemma}

Let us now define a space of analyzing vectors that allows us to 
treat also quasi-Banach spaces.
For $0<p\leq 1$ and for some submultiplicative weight function $w$ we define
\[
\BB^p_w \,:=\, \{g \in \D(K), V_g g \in W(L^\infty,L^p_w)\}.
\]
In the sequel we admit only those $p$ and $w$ such that $\BB^p_w \neq \{0\}$. 
Then it follows from the left and right translation invariance of $W(L^\infty,L^p_w)$ and the irreducibility
of $\pi$ that $\BB^p_w$ is dense in $\H$.
Now we are able to define the coorbit spaces.

\begin{definition}\label{def_coorbit} Let $Y$ be a solid quasi-Banach space of functions on $\G$ 
such that $W(L^\infty,Y)$ is
left and right translation invariant. Let $0<p\leq 1$ such that $Y$ has a $p$-norm and
put 
\begin{align}
\label{def_wweight}
w(x) \,&:=\, \max\{\on R_x|W(L^\infty,Y)\on, \Delta(x^{-1}) \on R_{x^{-1}}|W(L^\infty,Y)\on\},\\
\label{def_vweight}
v(x) \,& :=\, \max\{1, \on L_{x^{-1}} |W(L^\infty,Y)\on \}.
\end{align}
We assume that 
\begin{equation}\label{def_BY}
\BB(Y) \,:=\, \BB^p_w \cap \A_v
\end{equation}
 is non-trivial.
Then for $g \in \BB(Y) \setminus \{0\}$ 
the coorbit space is defined by
\[
\CC(Y) \,:=\, \Co W(L^\infty,Y) \,:=\, \{f \in (\H^1_v)^\urcorner, V_g f \in W(L^\infty,Y) \}
\]
with quasi-norm $\|f|\CC(Y)\| := \|V_g f|W(L^\infty,Y)\|$.
\end{definition}

Let us prove that the reproducing formula extends to $\CC(Y)$, and that $\CC(Y)$ is complete
and independent of the choice of $g \in \BB^p_w \setminus \{0\}$.

\begin{proposition}\label{prop_repform} 
Let $g \in \BB(Y)$ such that $\|Kg|\H\|=1$. A function 
$F \in W(L^\infty,Y)$ is of the form $V_g f$ for some $f \in \CC(Y)$ if and only if 
$F$ satisfies the reproducing formula $F=F*V_g g$.
\end{proposition}
\begin{Proof} If $f \in \CC(Y) \subset (\H^1_v)^\urcorner$ then the reproducing
formula $V_g f = V_g f * V_g g$ holds by the reproducing formula 
for $(\H^1_v)^\urcorner$, see Lemma \ref{lem_antidual}(c).

Conversely assume that $F=F*V_g g$ for some $F \in W(L^\infty,Y)$. 
By Lemma \ref{lem_Linfty_embed} $W(L^\infty,Y)$ is embedded into 
$L^{\infty}_{1/v}$.
Thus $F \in L^\infty_{1/v}$ and by Lemma \ref{lem_antidual}(d)
it holds $F= V_g f$ for some $f \in (\H^1_v)^\urcorner$, which is then automatically
contained in $\CC(Y)$ by assumption. 
\end{Proof}

\begin{Theorem}\label{thm_basic_prop}
\begin{itemize}
\item[(a)] $\CC(Y)$ is a quasi-Banach space.
\item[(b)] $\CC(Y)$ is independent of the choice of $g \in \BB(Y) \setminus \{0\}$.
\item[(c)] $\CC(Y)$ is independent of the reservoir $(\H^1_v)^\urcorner$
in the following sense: Assume that $\bS \subset \H^1_v$ is a locally convex
vector space, which is invariant under $\pi$,
and satisfies $\bS \cap \BB(Y) \neq \{0\}$. Denote by
$\bS^\urcorner$ its topological anti-dual. 
Then for a non-zero vector $g \in \BB(Y) \cap \bS$ it holds
$
\CC(Y) \,=\, \{f \in \bS^\urcorner: V_g f \in W(L^\infty,Y) \}.
$
\end{itemize}
\end{Theorem}
\begin{Proof} (a) 
Let $g \in \BB^p_w$ such that $\|Kg|\H\|=1$.
Assume $(f_n)_{n\in \N}$ to be a Cauchy sequence in $\CC(Y)$. This means
that $V_g f_n$ is a Cauchy sequence in $W(L^\infty,Y)$. By completeness of $W(L^\infty,Y)$ 
there exists a limit $F = \lim_{n \to \infty} V_g f_n$ in $W(L^\infty,Y)$. By Theorem 
\ref{thm_conv}(b) the definition of the weight $w$ implies that 
the operator $F \mapsto F * V_g g$ is continuous from $W(L^\infty,Y)$ into
$W(L^\infty,Y)$.  
Hence, we may interchange its application with taking limits, and together with the reproducing
formula (Proposition \ref{prop_repform}) this yields
\[
F \,=\, \lim_{n\to \infty} V_g f_n \,=\, \lim_{n\to \infty} V_g f_n * V_g g \,=\, F*V_g g.
\] 
Using Proposition \ref{prop_repform} once more we see that $F = V_g f$ for some
$f \in \CC(Y)$. Clearly, $f = \lim_{n \to \infty} f_n$ in $\CC(Y)$ and, hence, $\CC(Y)$ is complete.

(b) Let $g,g' \in \BB^p_w \setminus \{0\}$. Without loss of generality we may assume
that $g, g'$ are normalized, i.e., $\|Kg\|=\|Kg'\|=1$. 
Choose a vector $h \in \BB^p_w$ such that $Kh$ is
not orthogonal to $Kg$ and $Kg'$. It follows from the orthogonality relations
that
\[
0 \,\neq\, \langle Kg', Kh\rangle \langle Kh,Kg\rangle V_g g' \,=\, V_{g'} g' * V_h h * V_g g.
\] 
Since $V_g g^\nabla = V_g g$ and likewise for $h$ and $g'$, and since $w=w^*$ it follows from
Theorem \ref{thm_conv_Lp} that $V_g g' \in W(L^\infty,L^p_w)$. 
%
The inversion formula for $V_{g'}$ reads $g=\int_\G V_{g'}g(y)\pi(y) g' dy$, and one easily
deduces 
\begin{equation}\label{rep_form2}
V_g f = V_{g'} f * V_g g' \qquad \mbox{ for all } f \in (\H^1_v)^\urcorner.
\end{equation}
By the convolution relation in Theorem \ref{thm_conv}(b) we conclude 
$V_g f \in W(L^\infty,Y)$ if
$V_{g'} f \in W(L^\infty,Y)$. Exchanging the roles of $g$ and $g'$ shows the converse
implication.

(c) The $\pi$-invariance ensures that the voice transform 
$V_g f(x) = \langle f, \pi(x)g\rangle$ is well-defined for $f \in \bS^\urcorner$.
Moreover, it also implies that $\bS$ is dense in $\H$
by irreducibility of $\pi$. Thus from $\bS \subset \H^1_v$, it follows 
$(\H^1_v)^\urcorner \subset \bS^\urcorner$. Hence, the inclusion
$\CC(Y) \subset \{f \in \bS^\urcorner, V_g f \in W(L^\infty,Y)\}$ is clear.
Conversely, we have $W(L^\infty,Y) \subset L^\infty_{1/v}$ 
by Lemma \ref{lem_Linfty_embed} and thus
$V_g f \in W(L^\infty,Y)$ implies $V_g f \in L^\infty_{1/v}$. This is equivalent
to $f \in (\H^1_v)^\urcorner$, and the proof is completed. 
\end{Proof}

\begin{remark} The assumption that $W(L^\infty,Y)$ is left translation invariant
may even be weakened. Analyzing the previous proof it is enough to impose the 
following conditions: (a) $W(L^\infty,Y)$ is contained in a weighted $L^\infty_{1/v}$ space. (b) There exists a reservoir $\bS^\urcorner$ on which the wavelet transform 
is well-defined such that 
$\H_{1/v}^\infty:=\{f \in \bS^\urcorner, V_g f \in L^\infty_{1/v}\}$ 
satisfies the properties in Lemma \ref{lem_antidual}
(for some suitable $g$).
\end{remark}

Let us give a characterization of the space of analyzing vectors $\BB^p_w$.
For simplicity we restrict to the case that $\G$ is an IN group although
a similar (but more complicated) result can also be formulated 
in the general case.

\begin{Theorem}\label{BBpw_ident} 
Let $\G$ be an IN-group. 
Let $w$ be a 
submultiplicative weight and $0<p\leq 1$. Define 
$w^\bullet(x) := \max\{w(x),w(x^{-1})\} \geq 1$. Then it holds
\[
\BB^p_w = \BB^p_{w^\bullet} = \CC(L^p_{w^\bullet}).
\]
\end{Theorem}
\begin{Proof} 
Let $g \in \BB^p_w$. It follows from $V_g g= V_g g^\nabla$ and 
Lemma \ref{lem_IN_inv} that $V_g g \in W(L^\infty,L^p_{w^\bullet})$, i.e.,
$g \in \BB^p_{w^\bullet}$. 
Let $g'$ be another element of
$\BB^p_w=\BB^p_{w^\bullet}$. Then the previous proof
showed that $V_g g', V_{g'} g \in W(L^\infty,L^p_{w^\bullet})$ and thus
$g,g' \in \CC(L^p_{w^\bullet})$. 

Conversely, assume that 
$g \in \CC(L^p_{w^\bullet})$. Note that 
$\CC(L^p_{w^\bullet}) \subset \H = \D(K)$ (the latter by 
unimodularity of $\G$)
so that voice transforms are well-defined.
Let $g' \in \BB^p_{w^\bullet} \setminus \{0\}$.
Setting $f=g$ in (\ref{rep_form2}) shows 
$V_g g = V_{g'} g * V_g g' = V_{g'} g * (V_{g'} g)^\nabla$. Since both
$V_{g'} g$ and $(V_{g'} g)^\nabla$ are contained in 
$W(L^\infty,L^p_{w^\bullet})$ (Lemma \ref{lem_IN_inv}) it follows
from Theorem \ref{thm_conv_Lp} that 
$V_g g \in W(L^\infty,L^p_{w^\bullet})$, i.e., 
$g \in \BB^p_{w^\bullet}$.
\end{Proof}

The following theorem will be useful to prove a weak version of  a conjecture
in \cite[Conjecture 12]{GS}. 

\begin{Theorem}\label{thm_weakconjec} Let $\G$ be an IN group. 
Let $w$ be a submultiplicative weight function satisfying
$w = w^\vee$ and assume $0<p\leq 1$. If
$V_g f \in W(L^\infty,L^p_w)$ for $f,g \in \H$ then 
both $f$ and $g$ are contained in $\BB^p_w = \CC(L^p_w)$.
\end{Theorem}
\begin{Proof} Since $\G$ is unimodular, it holds $\H = \D(K)$.
It follows from (\ref{rep_form2}) that 
\[
V_g g = (V_g f)^\nabla * V_g f \quad  \mbox{and} \quad V_f f = V_g f * (V_g f)^\nabla.
\]
Since $w(x) = w(x^{-1})$ it follows from Lemma \ref{lem_IN_inv}
that also $(V_g f)^\nabla$ is contained in $W(L^\infty,L^p_w)$.
The convolution relation in Theorem
\ref{thm_conv_Lp} and Theorem \ref{BBpw_ident} finally yields the assertion.
\end{Proof}

\section{Discretizations}

Our next aim is to derive atomic decompositions for coorbit spaces.
The idea is to discretize the reproducing formula. For some suitable function $G$
 -- later on we will take $G= V_g g$ -- we denote the convolution operator
by
\[
T F \,:=\, T^G F \,=\, F*G \,=\, \int_\G F(y) L_y G\, dy. 
\]
For some BUPU $\Psi=(\psi_i)_{i\in I}$ with corresponding well-spread set $X=(x_i)_{i\in I}$ 
we define the approximation operator
\[
T_\Psi F \,:=\, T^\G_\Psi F \,:=\, \sum_{i\in I} \langle F, \psi_i\rangle L_{x_i} G.
\]
We first prove its boundedness.

\begin{proposition}\label{prop_an} Let $X=(x_i)_{i\in I}$ be some well-spread set
and $(\psi_i)_{i\in I}$ be a corresponding BUPU of size $U$. 
Then the operator
$F \mapsto (\langle F, \psi_i\rangle)_{i\in I}$ is continuous
from $W(L^\infty,Y)$ into $Y_d$, i.e.,
$
\|(\langle F, \psi_i\rangle)_{i\in I}|Y_d\| \leq C\|F|W(L^\infty,Y)\|.
$ 
\end{proposition}
\begin{Proof} Let $F \in W(L^\infty,Y)$ and $Q$ some compact 
neighborhood of $e$. Denoting $I_y:=\{i \in I: y \in x_i Q\}$ we obtain
\[
|\sum_{i\in I} \langle F, \psi_i\rangle \chi_{x_i Q}(y)|
\,\leq\, \sum_{i\in I_y} \langle |F|,\psi_i\rangle
\,\leq\, \langle L_y \chi_{Q^{-1} U}, |F|\rangle \,=\, |F| * \chi_{U^{-1}Q}(y).
\]
The function $\chi_{U^{-1}Q}$ is contained in $W(L^\infty,L^p_w)$ for any
$p$ and submultiplicative weight function $w$.
By solidity this yields
\begin{align}
\|(\langle F,\psi_i\rangle)_{i\in I}|Y_d\| 
\,&\leq \, \|F * \chi_{U^{-1}Q}|Y\| 
\,\leq\, \|F*\chi_{U^{-1}Q}|W(L^\infty,Y)\|\notag\\
&\leq\, C\|F|W(L^\infty,Y)\| \|\chi_{U^{-1}Q}|W(L^\infty,L^p_w)\|,\notag
\end{align}
where $p$ and $w$ are chosen according to Theorem \ref{thm_conv}(b).
\end{Proof} 

\begin{proposition}\label{prop_syn} Let $0<p\leq 1$ be such that $Y$ 
has a $p$-norm. Set $m(x) := \|A_x|W(M,Y)\|$ 
and $v(x) := \|L_{x^{-1}}|W(L^\infty,Y)\|$. 
Assume $G \in W(C_0,L^p_m) \cap W(C_0,L^1_v)$.
Then the mapping
\begin{equation}\label{def_lambda_series}
(\lambda_i)_{i\in I} \mapsto \sum_{i\in I} \lambda_i L_{x_i} G
\end{equation}
defines a bounded operator from $Y_d$ into $W(L^\infty,Y)$.
The sum always converges pointwise, and in the norm of $W(L^\infty,Y)$ if
the finite sequences are dense in $Y_d$.
\end{proposition}
\begin{Proof} Let $(\lambda_i)_{i\in I} \in Y_d$.
By Lemma 3.8(a) in \cite{FG2} the sequence $(L_{x_i}G(x))_{i\in I}$ is 
contained in $\ell^1_v$ and by the embedding $Y_d \subset \ell^\infty_{1/v}$
(Lemma \ref{lem_Linfty_embed}) the series in (\ref{def_lambda_series})
converges pointwise. 
We define the measure 
$\mu_\Lambda := \sum_{i\in I} \lambda_i \epsilon_{x_i}$ with 
$\epsilon_{x_i}$ being the Dirac measure at $x_i$. It follows from Theorem
\ref{thm_basic} (in particular (\ref{WYd_norm})) that $\mu_\Lambda$ is contained in $W(M,Y)$ with
$\|\mu_\Lambda|W(M,Y)\| \leq C \|(\lambda_i)_{i\in I}|Y_d\|$.
Observe that $\sum_{i \in I} \lambda_i L_{x_i} G = \mu_\Lambda * G$.
Using Theorem \ref{thm_conv}(a) we finally obtain
\begin{align}
&\|\sum_{i\in I} \lambda_i L_{x_i} G|W(L^\infty,Y)\| 
\,=\, \|\mu_\Lambda * G|W(L^\infty,Y)\|\notag\\ 
&\leq\, C\|\mu_\Lambda|W(M,Y)\| \, \|G|W(L^\infty,L^p_m)\|
\,\leq\, C' \|(\lambda_i)_{i\in I} |Y_d\| \, \|G|W(L^\infty,L^p_m)\|.\notag
\end{align}
This estimation also implies the norm convergence of the series in (\ref{def_lambda_series})
if the finite sequences are dense in $Y_d$.
\end{Proof}

\begin{corollary} Let $Y$, $p$, $m$ and $v$ be as in the previous proposition and assume
$G\in W(C_0,L^p_m) \cap W(C_0,L^1_v)$. 
Then the operator $T_\Psi = T^G_\Psi$ is bounded
from $W(L^\infty,Y)$ into $W(L^\infty,Y)$.
\end{corollary}

The following theorem will be the key for providing atomic decompositions.

\begin{Theorem}\label{thm_Tpsi_inv} Let $Y$ have 
a $p$-norm for $0<p\leq 1$ and let $w$ as in Definition
\ref{def_coorbit}. Suppose 
$G \in W(C_0, W(L^\infty,L^p_w)^\vee)^\vee\cap W(C_0,L^p_w)$. 
Further, assume $\Psi$ to be a BUPU of size $U$. Then it holds
\[
\on T^G - T^G_\Psi| W(L^\infty,Y) \on \,\leq\, C\|G_U^\#|W(C_0,L^p_w)\|.
\]
In particular, we have 
$\on T^G - T^G_\Psi|W(L^\infty,Y)\on \to  0$ if $U \to \{e\}$
by Lemma \ref{lem_osc}(b).
\end{Theorem}
\begin{Proof} For $F \in W(L^\infty,Y)$ we get
\begin{align}
|T F - T_\Psi F| \,=&\, \left| \sum_{i\in I} \int_\G F(y) \psi_i(y)
(L_y G - L_{x_i} G) dy \right| \notag\\
\,\leq&\, \sum_{i \in I} \int_\G |F(y)| \psi_i(y) |L_y G - L_{x_i} G| dy.\notag
\end{align}
Since $\supp \psi_i \in x_i U$ we obtain with Lemma \ref{lem_osc}(c) 
\[
|T F - T_\Psi F| \,\leq\, \sum_{i \in I} \int_\G |F(y)| \psi_i(y) L_y G_U^\# dy
\,=\, \int_\G |F(y)| L_y G_U^\# dy \,=\, |F| * G_U^\#.
\]
By Lemma \ref{lem_osc}(a) we have $G_U^\# \in W(C_0,L^p_w)$. Thus, 
the convolution relation in Theorem \ref{thm_conv}(b)
finally yields
\[
\|T F - T_\Psi F|W(L^\infty,Y)\|\,\leq\, C\|F|W(L^\infty,Y)\|\, \|G^\#_U|W(C_0,L^p_w)\|.
\]
This gives the estimate for the operator norm.
\end{Proof}

Before stating the general discretization theorem we need to introduce another
set of analyzing vectors. Let $v,w$ be the weight functions defined in 
(\ref{def_vweight}),(\ref{def_wweight}) and set
$
\tilde{w}(x) \,:=\, \max\{w(x),\on A_x|W(M,Y)\on\}.
$
Then we define
\begin{equation}\label{def_DD}
\DD(Y) := \A_v \cap 
\{g \in \H, V_g g \in W(L^\infty,L^p_{\tilde{w}}) \cap 
W(L^\infty,W(L^\infty,L^p_w)^\vee)^\vee\}.
\end{equation}
If $\G$ is an IN group and if $\|A_x|W(M,Y)\| \leq C w(x)$ 
then it holds $\DD(Y) = \BB(Y)$ by Remark \ref{rem_GIN}.

\begin{Theorem}\label{thm_atomic} Let $g \in \DD(Y) \setminus \{0\}$. 
Then there exists a compact 
neighborhood $U$ of $e$ such that for any $U$-dense well-spread set $X=(x_i)_{i\in I}$
the family $\{\pi(x_i) g\}_{i\in I}$ forms an atomic decomposition of $\CC(Y)$. 
This means that there exists a sequence $(\lambda_i)_{i\in I}$ of linear bounded functionals on 
$(\H^1_v)^\urcorner$ (not necessarily unique) such that
\begin{itemize}
\item[(a)] 
$
f = \sum_{i\in I} \lambda_i(f) \pi(x_i) g 
$
for all $f \in \CC(Y)$
with convergence in the weak-$*$ topology of $(\H^1_v)^\urcorner$, and in the quasi-norm
topology of $\CC(Y)$ provided the finite sequences are dense in $Y_d$;
\item[(b)] an element $f \in (\H^1_v)^\urcorner$ is contained in $\CC(Y)$ if and only if 
$(\lambda_i(f))_{i\in} \in Y_d$ and
\[
\|(\lambda_i(f))_{i\in I}|Y_d\| \,\asymp\, \|f|\CC(Y)\| \mbox{ for all } f\in \CC(Y).
\]
\end{itemize}
\end{Theorem}
\begin{Proof} Without loss of generality we may assume $\|Kg|\H\|=1$. 
Let $G:= V_g g$. Then 
$
G \in W(C_0, L^p_m) \cap W(C_0, L^1_v) \cap W(C_0,W(L^\infty,L^p_w)^\vee)^\vee
\cap W(C_0, L^p_w)
$
with $m(x):= \|A_x|W(M,Y)\|$ by assumption on $g$.
The restriction of 
the operator $T^G$ to the closed subspace $W(L^\infty,Y) * G$ is 
the identity by the reproducing
formula (Proposition \ref{prop_repform}). By 
Theorem \ref{thm_Tpsi_inv} we can find a
compact neighborhood $U$ of $e$ such that $\on T^G-T^G_\Psi|W(L^\infty,Y)\on < 1$ 
for 
any BUPU $\Psi$ of size $U$. 
Hence, $T_\Psi$ is invertible on $W(L^\infty,Y) * G$ by means of the von Neumann series.
Let $X=(x_i)_{i\in I}$ be the well-spread set corresponding to $\Psi$. 

If $f \in \CC(Y)$ then $V_g f \in W(L^\infty,Y)*G$ by the reproducing formula 
(Proposition \ref{prop_repform}) and
\[
V_g f \,=\, T^G_\Psi (T^G_\Psi)^{-1} V_g f \,=\, 
\sum_{i\in I} \langle (T_\Psi^G)^{-1} V_g f, \psi_i \rangle L_{x_i} G.
\]
We have $L_{x_i} G = L_{x_i} V_g g = V_g (\pi(x_i) g)$, see Lemma \ref{lem_antidual}(a). 
Since $V_g$ is an isometric isomorphism
between $\CC(Y)$ and $W(L^\infty,Y) * G$ we obtain
\[
f \,=\, \sum_{i\in I} \langle (T_\Psi^G)^{-1} V_g f,\psi_i \rangle \pi(x_i) g. 
\]
Set $\lambda_i(f) := \langle (T_\Psi^G)^{-1} V_g f, \psi_i\rangle$. Clearly, $\lambda_i$, $i\in I$, 
is a linear bounded functional on $\CC(Y)$ (and also on $(\H^1_v)^\urcorner$).
Since $(T_\Psi^G)^{-1} V_g f \in W(L^\infty,Y)*G \subset W(L^\infty,Y)$ Proposition
\ref{prop_an} yields
\begin{align}
& \|(\lambda_i(f))_{i\in I}|Y_d\| \,\leq\, C \|(T_\Psi^G)^{-1} V_g f|W(L^\infty,Y)\|\notag\\
&\leq\, C\on (T_\Psi^G)^{-1}|W(L^\infty,Y)*G\on \,\|V_g f|W(L^\infty,Y)\|
\,=\, C\on (T_\Psi^G)^{-1}\on \|f|\CC(Y)\|.\notag
\end{align}
Conversely, assume $(\lambda_i)_{i\in I} \in Y_d$.
We apply $V_g$ to the series $\sum_{i\in I} \lambda_i \pi(x_i) g$
to obtain (at least formally)
\[
F \,:=\, V_g \left( \sum_{i\in I} \lambda_i \pi(x_i) g\right) \,=\, 
\sum_{i\in I} \lambda_i L_{x_i} G.
\]
Since $Y_d \subset \ell^\infty_{1/v}$ and $G \in W(C_0, L^1_v)$ the series defining $F$ converges
pointwise and defines a function in $L^\infty_{1/v}$ by Lemma 3.8(a) in \cite{FG2}. 
The pointwise convergence of the partial
sums of $F$ implies the weak-$*$ convergence of $\sum_{i\in I} \lambda_i \pi(x_i) g = f$,
see Lemma \ref{lem_antidual}(a),(b).
Once $f$ is identified with an element of $(\H^1_v)^\urcorner$ it belongs
to $\CC(Y)$ by Proposition \ref{prop_syn}, i.e.,
\[
\|f|\CC(Y)\| \,=\, \|\sum_{i\in I} \lambda_i L_{x_i} G |W(L^\infty,Y)\| 
\,\leq\, C \|(\lambda_i)_{i\in I}|Y_d\|. 
\]
 Also the type of convergence follows from Proposition \ref{prop_syn}.
\end{Proof}

\begin{remark} By using similar techniques one can also
provide Banach frames of the form $\{\pi(x_i) g\}_{i\in I}$ for $\CC(Y)$,
compare also \cite{gro,hr2}.
\end{remark}

The following auxiliary result will be needed later on.

\begin{Theorem}\label{thm_sample} Let $g \in \BB(Y)$ and $X$ be some well-spread set. 
Then 
\[
\|(V_g f(x_i))_{i\in I}|Y_d\| \,\leq\, C \|f|\CC(Y)\|\qquad\mbox{ for all } 
f \in \CC(Y).
\]
\end{Theorem}
\begin{Proof} Let $f \in \CC(Y)$. 
The function $V_g f$ is continuous. Using Theorem \ref{thm_basic}(d)
we obtain
\begin{align}
\|(V_g f(x_i))_{i\in} |Y_d\| \,&\leq\, 
\|(\|V_g f \chi_{x_i U}\|_\infty)_{i\in I}|Y_d\| \,\leq\, 
C \|V_g f|W(L^\infty,Y)\|\notag\\
&=\, C\|f|\CC(Y)\|.\notag
\end{align}
\vskip-5ex
\end{Proof}

In certain situations one might be able to construct
certain expansions as in (\ref{expan_f}) below 
on the level of the Hilbert space $\H$. For instance, 
the construction of wavelet orthonormal bases of $L^2(\R^d)$ 
of this type (with $r=2^d-1$ tensor product wavelets) is well-known.  
Also tight frames of such kind have been constructed in time-frequency analysis
and wavelet analysis. The theorem below states that these expansions extend
automatically from $\H$ to general coorbit spaces under certain assumptions. 
Its proof is a modification of the one in
\cite{Gr3}.

\begin{Theorem}\label{thm_extend} Let $g_r, \gamma_r \in \BB(Y)$, $r=1,\hdots,n$, 
and $X=(x_i)_{i\in I}$ be a well-spread set such that
\begin{equation}\label{expan_f}
f\,=\, \sum_{r=1}^n \sum_{i\in I} 
\langle f, \pi(x_i) \gamma_r\rangle \pi(x_i) g_r \quad 
\mbox{for all } f \in \H.
\end{equation}
Then expansion (\ref{expan_f}) extends to all $f \in \CC(Y)$ with norm
convergence if the finite sequences are dense in $Y_d$ and with 
weak-$*$ convergence in general.
Moreover, $f \in (\H^1_v)^\urcorner$ is contained in $\CC(Y)$ if and only
if $(\langle f,\pi(x_i) \gamma_r\rangle)_{i\in I}$ is contained in $Y_d$ for each $r=1,\hdots,n$,
and 
\[
\|((\langle f, \pi(x_i) \gamma_r\rangle)_{i\in I})_{r=1}^n|\oplus_{r=1}^n Y_d\|
\asymp \|f|\CC(Y)\| \quad \mbox{ for all } f\in \CC(Y).
\]
\end{Theorem}
\begin{Proof} We prove the theorem for the case that the finite sequences
are dense in $Y_d$. The general case requires some slight changes.

Let $g \in \DD(Y)$ and $Z=(z_j)_{j\in J}$ be a well-spread set as in Theorem \ref{thm_atomic},
i.e., such that $\{\pi(z_j)g\}_{j\in J}$ forms an atomic decomposition of $\CC(Y)$.
Let $f \in \CC(Y)$. Then we have the decomposition $f= \sum_{j\in J} \lambda_j(f) \pi(z_j)g$
with $(\lambda_j(f))_{j\in J} \in Y_d(Z)$. Since the finite sequences are dense in $Y_d$ there
exists a finite set $N \subset J$ such that $\|\Lambda_N|Y_d(Z)\| < \epsilon$ where 
$(\Lambda_N)_j = \lambda_j$ if $j \notin N$ and $(\Lambda_N)_j = 0$ otherwise.
The element $f_N \,=\, \sum_{j\in N} \lambda_j \pi(z_j) g$ is contained in $\H \cap \CC(Y)$
and satisfies 
\[
\|f-f_N|\CC(Y)\| \,\leq\, C \|\Lambda_N|Y_d(Z)\| \,\leq\, C\epsilon.
\]
By assumption $f_N$ has the expansion
\[
f_N \,=\, \sum_{r=1}^n \sum_{i\in I} \langle f_N, \pi(x_i) \gamma_r \rangle \pi(x_i)g_r.
\] 
Theorem \ref{thm_sample} asserts that the sequence 
$(\langle f_N, \pi(x_i) \gamma_r\rangle)_{i\in I} = (V_{\gamma_r} f_N(x_i))_{i\in I}$ is 
contained in $Y_d(X)$ for all $r=1,\hdots,n$. As above there exist finite 
sets $N_r$, $r=1,\hdots,n$, such that $\|\kappa^{(r)}|Y_d(X)\| \leq \epsilon$
where $\kappa^{(r)}_i = \langle f_N, \pi(x_i) \gamma_r \rangle$ if $i \notin N_r$ and 
$\kappa^{(r)}_i = 0$ otherwise. Then
$
f_1 := \sum_{r=1}^n \sum_{i \in N_r} \langle f_N, \pi(x_i) \gamma_r\rangle \pi(x_i) g_r
$
satisfies 
\[
\|f_N-f_1|\CC(Y)\| \leq C' \sum_{r=1}^n \|\kappa^{(r)}|Y_d(X)\| \leq Cr\epsilon.
\]
It follows from Proposition \ref{prop_syn} that 
\begin{align}
&\|\sum_{r=1}^n \sum_{i\in I} \lambda_{i,r} \pi(x_i) g_r|\CC(Y)\|^p
\,\leq\, \sum_{r=1}^n C_r \|\sum_{i \in I} \lambda_{i,r} L_{x_i} V_{g_r} g_r|W(L^\infty,Y)\|^p \notag\\
&\leq\, C \sum_{r=1}^n \|(\lambda_{i,r})_{i\in I}|Y_d(X)\|^p
\,\leq\, C' \|((\lambda_{i,r})_{i\in I})_{r=1}^n|\oplus_{r=1}^n Y_d(X)\|^p \notag
\end{align}
for all 
$(\lambda_{i,r}) \in \oplus_{r=1}^n Y_d(X)$.
With $f_2 := \sum_{r=1}^n \sum_{i\in N_r} \langle f, \pi(x_i) \gamma_r\rangle \pi(x_i) g_r$
this yields using Theorem \ref{thm_sample}
\begin{align}
\|f_1 - f_2|\CC(Y)\| \,&\leq\, 
C\|((\langle f_N-f, \pi(x_i)\gamma_r\rangle)_{i\in I})_{r=1}^n|\oplus_{r=1}^n Y_d\|\notag\\
&\leq\, C''\|f_N-f|\CC(Y)\| \,\leq\, C''\epsilon. \notag
\end{align}
Altogether we get
\begin{align}
\|f-f_2|\CC(Y)\| \,&\leq\, \|f-f_N|\CC(Y)\| + \|f_N-f_1|\CC(Y)\| + \|f_1 -f_2|\CC(Y)\|\notag\\
&\leq\, (C+C'r+C'')\epsilon.\notag
\end{align}
Since $\epsilon$ can be chosen arbitrarily small we deduce that
\[
f = \sum_{r=1}^n \sum_{i\in I} \langle f,\pi(x_i) \gamma_r\rangle \pi(x_i) g_r
\]
for all $f \in \CC(Y)$ with quasi-norm convergence. Moreover, 
it follows from Theorem \ref{thm_sample}
and Proposition \ref{prop_syn} that
\begin{align}
\|f|\CC(Y)\| \,&=\, \|\sum_{r=1}^n \sum_{i\in I} \langle f, \pi(x_i) \gamma_r\rangle \pi(x_i) g_r|\CC(Y)\| \notag\\
\,&\leq\,  C (\sum_{r=1}^n \|(\langle f,\pi(x_i) \gamma_r)_{i\in I}|Y_d(X)\|^p)^{1/p}\notag\\
&\leq\, C (\sum_{r=1}^n \|V_{\gamma_r} f|W(L^\infty,Y)\|^p)^{1/p}
\,\leq \, C' \|f|\CC(Y)\|.\notag
\end{align}
This concludes the proof.
\end{Proof}

\section{Characterizations of $\CC(Y)$ via $Y$}

The original definition of the coorbit spaces by Feichtinger and Gr\"ochenig
involves $Y$ rather than $W(L^\infty,Y)$. It is interesting to investigate what
happens if we replace $W(L^\infty,Y)$ by $Y$ in our more general case. 
In order to 
distinguish clearly between the two spaces let us
denote $\Co Y = \{f \in (\H^1_v)^\urcorner, V_g f \in Y\}$ with natural norm
$\|f|\Co Y\| = \|V_g f|Y\|$, and $\Co W(L^\infty,Y) = \CC(Y)$ as usual.
It was already proven in \cite{FG3} that in the classical 
Banach space case both
spaces coincide:

\begin{Theorem} (Theorem 8.3 in \cite{FG3}) Let $Y$ be a solid Banach space of 
functions on $\G$ that is left and right translation invariant and continuously
embedded into $L^1_{loc}(\G)$. Then it holds 
$\Co Y = \Co W(L^\infty,Y)$
with equivalent norms. 
\end{Theorem}

In the general case of quasi-Banach spaces at least the inclusion
$\Co W(L^\infty,Y) \subset \Co Y$ holds 
since $W(L^\infty,Y) \subset Y$.
However, it seems doubtful that we can state results on the converse
inclusion in the general abstract case. Indeed, in order 
to come up with an abstract result
we would need a convolution relation of the type 
$Y*W(L^\infty,L^p_w) \hookrightarrow W(L^\infty,Y)$
(see \cite{FG2} for the Banach space case).
However, such a relation does not
hold for general quasi-Banach spaces (consider $Y=L^p(\G)$, $0<p<1$, for a 
non-discrete group $\G$).

Moreover, it is even not clear whether $\Co Y$
is a complete space. Indeed, in the proof of completeness of $\CC(Y)$ (and of
$\Co Y$ in the case of Banach spaces $Y$) one makes heavy use of 
the reproducing
formula together with the fact that convolution from the right with a suitable $G$
is a bounded operator from $W(L^\infty,Y)$ into itself. 

In special cases, however, one might be able to prove that 
$\|V_g f|W(L^\infty,Y)\| \leq C \|V_g f|Y\|$ for a very specific choice of $g$,
by using methods that are not available in the abstract setting (like analyticity 
properties for instance), see e.g. Section 7. 
Then one may extend this inequality 
to more general analyzing vectors $g$ 
as shown by the next result.

\begin{Theorem}\label{thm_charY} Let $Y$ be a right translation 
solid $p$-normed quasi-Banach space of 
functions on $\G$ such that $W(L^\infty,Y)$ is left translation invariant.
Set $\nu(x) := \on R_x|Y\on$. 
Let $\bS \subset \H^1_v$ be a $\pi$-invariant 
locally convex vector space (with 
$v(x):=\max\{1,\on L_{x^{-1}}|W(L^\infty,Y)\on\}$, see (\ref{def_vweight})). 
Assume that there exists a non-zero vector 
$g_0 \in \BB(Y) \cap \CC(L^p_\nu) \cap \bS$ and a constant $C>0$ such that
\[
\|V_{g_0} f |W(L^\infty,Y)\| \,\leq\, C \|V_{g_0} f|Y\|
\]
for all $f \in \bS^\urcorner \supset (\H^1_v)^\urcorner$ (with the 
understanding that $V_{g_0} f \in W(L^\infty,Y)$
if $V_{g_0} f \in Y$).
Let 
$g \in \BB(Y) \cap \DD(L^p_\nu) \cap \bS \setminus \{0\}$. Then it holds
\[
\|V_g f|W(L^\infty,Y)\| \,\asymp \,  \|V_g f|Y\|
\]
for all $f \in \bS^\urcorner \supset (\H^1_v)^\urcorner$ and
$
\Co W(L^\infty,Y) = \Co Y
= \{f \in \bS^\urcorner, V_g f \in Y\}.
$ 
In particular, $\Co Y$ is complete with the quasi-norm $\|V_g f|Y\|$ 
and independent of the choice of 
$g \in \BB(Y) \cap \DD(L^p_\nu) \cap \H^1_v \setminus \{0\}$.
\end{Theorem}
\begin{Proof} Since $\CC(Y)$ is independent of the 
choice of $g \in \BB(Y) \setminus \{0\}$ and
of the reservoir $\bS$ (Theorem \ref{thm_basic_prop}) we have 
\[
\|V_g f|W(L^\infty,Y)\| \leq C \|V_{g_0} f|W(L^\infty,Y)\|.
\]
for all $f \in \bS^\urcorner$.
Thus, it remains to prove  that
$\|V_{g_0} f|Y\| \leq C \|V_g f|Y\|$ for all $f \in \bS^\urcorner$.
By the assumptions on $g$ it follows from 
Theorem \ref{thm_atomic} that $g_0$ has a decomposition
\[
g_0 \,=\, \sum_{i\in I} \lambda_i(g_0) \pi(x_i) g
\]
with $(\lambda_i(g_0))_{i\in I} \in \ell^p_{\nu} = (L^p_{\nu})_d$
and  
$\|(\lambda_i(g_0))_{i\in I}|\ell^p_{\nu}\| \asymp \|g_0|\CC(L^p_\nu)\|$.
Hence, we obtain
\[
V_{g_0} f(x) \,=\, \langle f, \pi(x) g_0 \rangle 
\,=\, \langle f, \pi(x) \sum_{i\in I} \lambda_i(g_0) \pi(x_i) g\rangle 
\,=\, \sum_{i\in I} \ol{\lambda_i(g_0)} R_{x_i} V_g f(x).
\]
This yields
\begin{align}
\|V_{g_0} f|Y\|^p 
\,&=\, \|\sum_{i\in I} \ol{\lambda_i(g_0)} R_{x_i} V_g f |Y\|^p 
\,\leq\, \sum_{i\in I} |\lambda_i(g_0)|^p \on R_{x_i}|Y\on^p \|V_g f|Y\|^p
\notag\\
&\leq\, C \|g_0|\CC(L^p_v)\|^p\, \|V_g f|Y\|^p\notag
\end{align}
for all $f \in \bS^\urcorner$. The reverse inequality 
$\|V_g f|Y\| \leq \|V_g f|W(L^\infty,Y)\|$ is clear.
\end{Proof}

\section{Nonlinear Approximation}

Let us now discuss non-linear approximation. 
Let $(x_i)_{i\in I}$ be some well-spread set and $g$ 
such that $\{\pi(x_i)g\}_{i\in I}$ forms an atomic decomposition of the coorbit
space we want to consider.
We denote by
\[
\sigma_n(f,\CC(Y)) \,:=\, \inf_{N \subset I, \#N \leq n} \|f - \sum_{i \in N} \lambda_i \pi(x_i) g|\CC(Y)\|
\]
the error of best $n$-term approximation in $\CC(Y)$. 
Hereby, the infimum is also taken over all possible
choices of coefficients $\lambda_i$. 
Our task is to find a class of elements for which this error has a certain decay
when $n$ tends to $\infty$.

To this end we consider coorbits with respect to Lorentz spaces. 
For some measurable function $F$ on $\G$ let $\lambda_F(s) = |\{x: |F(x)|>s\}|$ denote its distribution
function and $F^*(t) = \inf \{s: \lambda_F(s) \leq t\}$ 
its non-increasing rearrangement. 
For $0<p,q <\infty$ we let
\begin{align}
\|F\|_{p,q}^* \,&:=\, \left(\frac{q}{p}\int_0^\infty F^*(t)^q t^{q/p} \frac{dt}{t}\right)^{1/q}\notag\\
\mbox{and}\qquad 
\label{star_qnorm}
\|F\|_{p,\infty}^* \,&:=\, \sup_{t>0} t^{1/p} F^*(t).
\end{align}
The Lorentz-space $L(p,q)$ is defined as the collection of all $F$ such that the quantity above is finite.
It is well-known that $L(p,p) = L^p$. Most interesting for our purposes is $L(p,\infty)$, 
which is also called weak-$L^p$. It holds $L^p \subset L(p,\infty)$.
Sometimes it is useful to introduce another quasi-norm on $L(p,q)$. 
For some $r$ satisfying $0<r\leq 1$, $r < p$ and $r \leq q$  we let
\begin{equation}\label{def_Fss}
F^{**}(t,r) \,:=\, \sup_{t<\mu(E)<\infty} \left(\frac{1}{\mu(E)}\int_E |F(x)|^r dx\right)^{1/r},
\end{equation}
the supremum taken over measurable subsets of $\G$ with the stated property. We define
\[
\|F\|_{p,q}^{(r)} \,:=\, \left(\frac{q}{p}\int_0^\infty F^{**}(t,r)^q t^{q/p} \frac{dt}{t}\right)^{1/q}
\]
with modification for $q = \infty$ as in (\ref{star_qnorm}). 
It can be easily seen from (\ref{def_Fss}) that $\|\cdot\|_{p,q}^{(r)}$ is a quasi-norm.
If $r =1$ (implying $p>1$) and $q\geq 1$ then it is even a norm, and 
if $r<1$ and $q\geq 1$ then 
$\|\cdot\|_{p,q}^{(r)}$ is an $r$-norm.
Furthermore, if it can be shown \cite{Hunt} that
\[
\|F\|^*_{p,q} \,\leq\, \|F\|_{p,q}^{(r)} \,\leq\, \left(\frac{p}{p-r}\right)^{1/r} \|F\|^*_{p,q}.
\] 
In particular, also $\|F\|^*_{p,q}$ is a quasi-norm. 
Moreover, if $p>1$ and $q\geq 1$ then $L(p,q)$ with
the equivalent norm $\|\cdot\|_{p,q}^{(1)}$ is a Banach space. 
In the general case, $L(p,q)$ is only a quasi-Banach space.
Indeed, it is known that $L(1,q)$ is not normable 
for $q>1$ (except for the trivial case that the underlying
group is finite).

By the properties of the Haar-measure it is easily seen that all spaces $L(p,q)$ are 
left and 
right translation invariant. Thus, if $m$ is a moderate function then also 
$$
L_m(p,q) \,=\, \{F \mbox{ measurable}, Fm \in L(p,q) \}
$$
with the quasi-norm $\|F|L_m(p,q)\| := \|Fm\|_{p,q}^*$ is left and right translation invariant. 
In particular, the Wiener amalgam spaces $W(L^\infty,L_m(p,q))$ are well-defined. Further, 
if $\BB(L_m(p,q))$, see (\ref{def_BY}), is non-trivial then also the coorbit space 
$\CC(L_m(p,q))$ is well-defined. 

It is not difficult to see that the sequence space $(L_m(p,q))_d(X)$ associated to 
a well-spread set $X=(x_i)_{i\in I}$ coincides with a Lorentz space $\ell_m(p,q)$ 
on the index set $I$.
In particular, an equivalent quasi-norm on $(L_m(p,\infty))_d(X)$ is given by
\begin{equation}\label{Lorentz_sequence}
\|(\lambda_i)_{i\in I}\|_{p,\infty}^* \,=\, \sup_{n \in \N} n^{1/p}(\lambda m)^*(n)
\end{equation}
where $(\lambda m)^*$ denotes the non-increasing rearrangement of the sequence
$(\lambda_i m(x_i))_{i\in I}$.

\begin{Theorem}\label{thm_nonlinapprox} 
Let $m$ be some $w$-moderate weight function on $\G$, 
let $0 < p < q \leq \infty$
and define $\alpha = 1/p-1/q > 0$.
Let $(x_i)_{i\in I}$ be some well-spread set and 
$g \in \DD(L_m(p,\infty)) \subset \DD(L^q_m)$
such that $\{\pi(x_i)g\}_{i\in I}$ forms an atomic decomposition
simultaneously of $\CC(L_m(p,\infty))$ and $\CC(L^q_m)$ 
(according to Theorem \ref{thm_atomic}).
Then for all $f\in \CC(L_m(p,\infty))$ it holds
\begin{equation}\label{best_nterm_decay}
\sigma_n(f,\CC(L^q_m)) \,\leq\, C \|f|\CC(L_m(p,\infty))\|\, n^{-\alpha}.
\end{equation}
\end{Theorem} 
\begin{Proof}
Let $f = \sum_{i\in I} \lambda_i(f) \pi(x_i)g$ be an expansion
of $f \in \CC(L_m(p,\infty))$ in terms of the atomic decomposition.
By Theorem \ref{thm_atomic} it holds 
$(\lambda m)^*_k \leq C\|f|\CC(L_m(p,\infty))\| k^{-1/p}$. Let
$\tau: \N\to I$ be a bijection that realizes the non-increasing rearrangement, 
i.e., $\lambda_{\tau(k)} m(x_{\tau(k)}) = (\lambda m)^*_k$.
Moreover, $(L^q_m)_d = \ell^q_m(I)$, and 
$\|(\lambda_i(f))_{i\in I}|\ell^q_m(I)\|$ forms
an equivalent norm on $\CC(L^q_v)$ once again by Theorem \ref{thm_atomic}. 
We obtain
\begin{align}
&\sigma_n(f,\CC(L^q_m)) \,\leq\, 
\|f- \sum_{k=1}^n \lambda_{\tau(s)} \pi(x_{\tau(s)})g|\CC(L^q_m)\|\notag\\
&=\, \|\sum_{k={n+1}}^\infty \lambda_{\tau(k)} \pi(x_{\tau(k)}) g |\CC(L^q_m)\|
\,\leq\, C \left(\sum_{k=n+1}^\infty ((\lambda m)^*_k)^q\right)^{1/q}\notag\\
&\leq\, C \|f|\CC(L_m(p,\infty))\| \left(\sum_{k=n+1} k^{-q/p}\right)^{1/q}
\,\leq\, C \|f|\CC(L_m(p,\infty))\| (n^{-q/p+1})^{1/q}\notag\\
&=\,  C \|f|\CC(L_m(p,\infty))\| n^{-\alpha}.\notag
\end{align}
This completes the proof.
\end{Proof}
\begin{remark}
\begin{itemize}
\item[(a)] The obvious embedding $\CC(L^p_m) \,\subset\, \CC(L_m(p,\infty))$
implies that we also have $\sigma_n(f,\CC(L^q_m)) \leq Cn^{1/q-1/p}$
for all $f \in \CC(L^p_m)$ if $p<q$. However, in this situation one may even
prove a slightly faster decay of $\sigma_n(f,\CC(L^p_m))$, i.e, 
$o(n^{1/q-1/p})$ instead of $\O(n^{1/q-1/p})$, with methods similar as in 
\cite{GroS} for instance.
\item[(b)] In order to have a very fast decay of $\sigma_n(f,\CC(L^q_m))$ one 
obviously has to take $p$ very small in the Theorem above, in particular, 
$p\leq 1$.
Clearly, $\CC(L_m(p,\infty))$ is no longer a Banach space in this case, but only a 
quasi-Banach space. So it is very natural to treat also the case of 
quasi-Banach spaces when dealing with problems in non-linear approximation. This
was actually one of the motivations for this paper.
\item[(c)] The most interesting case appears when taking $q=2$ and 
$m\equiv 1$ since $\CC(L^2) = \H$. So for all $f \in \CC(L(p,\infty))$, $p<2$, we have
\[
\sigma_n(f,\H) \,\leq\, C \|f|\CC(L(p,\infty))\| \, n^{-1/p+1/2}.
\]
\end{itemize}
\end{remark}

\section{Modulation spaces}

Let $\HH_d := \R^d \times \R^d \times \TT$ denote the (reduced) Heisenberg group with group law
\[
(x,\omega,\tau) (x',\omega',\tau') ~=~ 
(x+x',\omega + \omega', \tau \tau' e^{\pi i(x'\cdot \omega - x\cdot \omega')}).
\]
It is unimodular and has Haar measure 
\[
\int_{\HH_d} f(h) dh = \int_{\R^d} \int_{\R^d} \int_0^1 f(x,\omega,e^{2\pi i t}) dt d\omega dx.
\]
We denote by 
\[
T_x f(t) \,:=\, f(t-x), \quad \mbox{and}\quad M_\omega f(t) \,=\, e^{2\pi i \omega \cdot t} f(t), \quad
x,\omega, t \in \R^d,
\]
the translation and modulation operator on $L^2(\R^d)$. 
Then the Schr\"odinger representation $\rho$ is defined by
\[
\rho(x,\omega,\tau) \,:=\, \tau e^{\pi i x\cdot \omega} T_x M_\omega
\,=\, \tau e^{-\pi i x \cdot \omega} M_\omega T_x.
\]
It is well-known that this is an irreducible unitary and square-integrable representation of $\HH_d$.
The corresponding voice transform is essentially the short time Fourier 
transform:
\begin{align}
&V_g f(x,\omega,\tau) \,=\, \langle f, \rho(x,\omega,\tau) g\rangle_{L^2(\R^d)} 
\,=\, \overline{\tau} 
\int_{\R^d} f(t) \overline{e^{-\pi i x\cdot \omega} M_\omega T_x g(t)} dt\notag\\
\label{Schroed_STFT}
&=\, \overline{\tau} e^{\pi i x\cdot \omega} 
\int_{\R^d} f(t) \overline{g(t-x)} e^{-2\pi i t \cdot \omega} dt 
\,=\, \overline{\tau} e^{\pi i x\cdot \omega} \STFT_g f(x,\omega).
\end{align}

Let us now introduce the modulation spaces on $\R^d$. We consider nonnegative
continuous weight functions $m$ on $\R^d \times \R^d$ that satisfy
\[
m(x+y,\omega + \xi) \,\leq\, C (1 + |x|^2 + |\omega|^2)^{a/2} m(y,\xi), \quad (x,\omega),(y,\xi) \in \R^d \times \R^d.
\]
for some constants $C>0,a \geq 0$. 
This means that $m$ is a $w$-moderate function for $w(x,\omega) = (1+|x|^2+|\omega|^2)^{a/2}$,
see also \cite[Chapter 11.1]{Gr}.
Additionally, we require $m$ to be symmetric, i.e.,
$m(-x,-\omega) = m(x,\omega)$.
A typical choice is 
$m_s(x,\omega) = (1+|\omega|)^{s}, s \in \R$.
For $0<p,q\leq \infty$ and $m$ as above we introduce 
$
L^{p,q}_m:= L^{p,q}_m(\R^{2d}) := \{F \mbox{ measurable}, \|F|L^{p,q}_m\| < \infty\}$
with quasi-norm
\[
\|F|L^{p,q}_m\| \,:=\, \left(\int_{\R^d}\left(\int_{\R^d} |F(x,\omega)|^p m(x,\omega)^p dx\right)^{q/p}d\omega\right)^{1/q}. 
\]
This expression is an $r$-norm with $r:=\min\{1,p,q\}$.

Let $g$ be some non-zero Schwartz function on $\R^d$. 
The short time Fourier transform $\STFT_g$ extends to the space $\cS'(\R^d)$ of tempered distributions
in a natural way. 
Given $0 <  p,q \leq \infty$ and $m$ as above the modulation space is defined
as 
\[
M^{p,q}_m \,:=\, \{f \in \cS'(\R^d), \|\STFT_g f|L^{p,q}_m\| < \infty\}
\]
with quasi-norm $\|f|M^{p,q}_m\| = \|\STFT_g f|L^{p,q}_m\|$.
Since by (\ref{Schroed_STFT}) $|V_gf(x,\omega,\tau)| = |\STFT_g f(x,\omega)|$ 
we can identify the modulation
spaces with coorbit spaces,
\[
M^{p,q}_m(\R^d) \,=\, \Co L^{p,q}_m(\HH_d) \,=\, \{f \in \cS', V_g f \in L^{p,q}_m\},
\]
where $m$ and $L^{p,q}_m$ are extended to $\HH_d$ in an obvious way, 
e.g.  $m(x,\omega,\tau) = m(x,\omega)$.
However, at the moment we do not know yet whether $\Co L^{p,q}_m$ coincides with 
\[
\CC(L^{p,q}_m) = \{f \in \cS'(\R^d), V_g f \in W(L^\infty,L^{p,q}_m)\}
\] 
if $p<1$ or $q<1$. It is even not clear yet whether $M^{p,q}_m$ is complete.
We will use Theorem \ref{thm_charY} and a result from \cite{GS} to clarify this problem.
Let us first investigate the spaces $\BB(L^{p,q}_m)$ and $\DD(L^{p,q}_m)$, see 
Definition \ref{def_coorbit} and (\ref{def_DD}).
One easily shows \cite[Lemma 4.7.1]{hr} (using the symmetry of $w$) that
$\on L_{(x,\omega,\tau)}|L^{p,q}_m\on \leq w(x,\omega)$ and 
$\on R_{(x,\omega,\tau)}|L^{p,q}_m\on \,\leq\, w(x,\omega)$.
Since $\HH_d$ is an IN group Lemma \ref{lem_IN_trans} yields
\[
\on R_{(x,\omega,\tau)}|W(L^\infty,L^{p,q}_m)\on \,\leq\, w(x,\omega)
\quad\mbox{and}\quad 
\on A_{(x,\omega,\tau)}|W(M,L^{p,q}_m)\on \,\leq\, w(x,\omega).
\]
Further, $\on L_{(x,\omega,\tau)}|W(L^\infty,L^{p,q}_m)\on \leq w(x,\omega)$ by
Lemma \ref{lem_left_trans}. 
Thus, using Remark \ref{rem_GIN} we conclude
\[
\BB^r_w \subset \BB(L^{p,q}_m) 
\quad \mbox{and}\quad \BB^r_w \subset \DD(L^{p,q}_m)\quad \mbox{with } 
r:= \min\{1,p,q\}.
\]
Moreover, Theorem \ref{BBpw_ident} yields 
\[
\BB^r_w \,=\, \CC(L^r_w).
\]

Let $g_0(t) = e^{-\pi |t|^2}$ be a Gaussian. 
Using the relation of $\STFT_{g_0}$ to the Bargmann transform Galperin and Samarah
proved that 
\[
\|V_{g_0} f|W(L^\infty,L^{p,q}_m)\| \,\leq\, C \|V_{g_0} f|L^{p,q}_m\|
\]
for all $f \in M^{p,q}_m$ \cite[Lemma 3.2]{GS}. 
Moreover, it is clear that \linebreak
$g_0 \in  \CC(L^r_w) \subset \BB(L^{p,q}_m) \cap \CC(L^p_\nu)$ and 
$\cS \subset \CC(L^r_w) \subset\BB(L^{p,q}_m) \cap \DD(L^p_\nu)$ with 
$\nu(x,\omega) := \on R_{(x,\omega)}|W(L^\infty,L^{p,q}_m)\on$.
Thus, by Theorem \ref{thm_charY}
\[
\CC(L^{p,q}_m) \,=\, M^{p,q}_m,
\]
and the latter is complete. It seems that the completeness of $M^{p,q}_m$ 
for $p<1$ or $q<1$ was not stated in \cite{GS} or elsewhere in the literature 
although its proof is somehow hidden in \cite{GS}.

The abstract discretization Theorem \ref{thm_atomic} 
yields the following result for atomic decompositions
of modulation spaces.

\begin{Theorem}\label{thm_Mod_atomic} 
Let $0<p_0 \leq 1$ and $w$ be some symmetric 
submultiplicative weight function
on $\R^d \times \R^d$ with polynomial growth. Assume $g \in M^{p_0}_w$. 
Then there exist constants $a,b>0$ such that
\[
\{M_{bj} T_{ak} g: k,j \in \Z^{d}\}
\]
forms an atomic decomposition for all modulation spaces $M^{p,q}_m$
with $p_0 \leq p,q \leq \infty$ and $m$ being a $w$-moderate weight. This means that there
exist functionals $\lambda_{k,j}$, $k,j \in \Z^d$ on $M^\infty_{1/w}$
$(\subset \cS')$ such that
\begin{itemize}
\item[(a)] any $f \in M^{p,q}_m$ has the series expansion 
$
f \,=\, \sum_{k,j \in \Z^d} \lambda_{k,j}(f) M_{bj} T_{ak} g\,;
$ 
\item[(b)] a distribution $f \in M^\infty_{1/w}$ belongs to $M^{p,q}_m$ if and
only if 
$(\lambda_{k,j}(f))_{k,j \in \Z^d}$
belongs to $\ell^{p,q}_m(\Z^{2d})$, and we have the quasi-norm equivalence
\begin{align}
\|f|M^{p,q}_m\| \,&\asymp\, 
\left(\sum_{j \in \Z^d} \left(\sum_{k \in \Z^d} |\lambda_{k,j}(f)|^p m(ak,bj)^p\right)^{q/p}\right)^{1/q} \notag\\
&=:\, \|(\lambda_{k,j}(f))|\ell^{p,q}_m(\Z^{2d})\|.\notag
\end{align}
\end{itemize}
\end{Theorem}  
  
We remark that the abstract Theorem \ref{thm_atomic} allows to extend the previous 
result also to irregular Gabor frames on $M^{p,q}_m$. 
 
Theorem \ref{thm_Mod_atomic} 
indicates that the modulation spaces $M^{p_0}_w$ with $0<p_0 \leq 1$
are the correct window classes for time-frequency analysis on $M^{p,q}_m$. This
was already conjectured in \cite{GS}. Galperin and Samarah also conjectured
that whenever $\STFT_g f \in L^p_v$ then $f \in M^p_v$ and $g \in M^p_v$ 
\cite[Conjecture 12]{GS}. 
Theorem \ref{thm_weakconjec} leads to a weak version of this conjecture.

\begin{Theorem} Let $f,g \in L^2(\R^d)$ and $0<p\leq 1$. Assume $v$ to be a symmetric
submultiplicative weight function. 
If $\STFT_g f \in W(L^\infty,L^p_v)$ 
then $g \in M^p_v$ and $f \in M^p_v$.
\end{Theorem}
The remaining question is whether $V_g f \in L^p_v$ already implies 
$V_g f \in W(L^\infty,L^p_v)$. 

Let us also apply Theorem \ref{thm_extend} to our situation.

\begin{Theorem} Let $g \in \cS(\R^d)$ and $a,b > 0$ such that 
\begin{equation}\label{Gab_frame}
\{M_{bj} T_{ak} g:\, j,k \in \Z^d\}
\end{equation}
forms a Gabor frame for $L^2(\R^d)$. Then its canonical dual $\gamma$ is also
contained in $\cS(\R^d)$, and any $f \in M^{p,q}_m$, $0 < p,q \leq \infty$, has
a decomposition
\[
f \,=\, \sum_{j,k \in \Z^d} \langle f, M_{bj} T_{ak} g\rangle M_{bj} T_{ak} \gamma
\]
with
$
\|f|M^{p,q}_m\| \,\asymp\,
\|(\langle f, M_{bj} T_{ak} g\rangle)_{j,k \in \Z^d}|\ell^{p,q}_m(\Z^{2d})\|.
$
\end{Theorem}
\begin{Proof} Since (\ref{Gab_frame}) forms a Gabor frame with dual window $\gamma$ 
any $f \in L^2(\R^d)$ has a decomposition
\[
f \,=\, \sum_{j,k \in \Z^d} \langle f, M_{bj} T_{ak} g\rangle M_{bj} T_{ak} \gamma.
\]
It was shown in \cite[Corollary 13.5.4]{Gr} that also the dual window 
$\gamma$ is contained in $\cS(\R^d)$. Since $\cS(\R^d) \subset M^{p,q}_m$ for all
$0<p,q\leq \infty$ and all $w$-moderate weights $m$ with $w$ having 
polynomial growth we have 
$g,\gamma \in \BB(L^{p,q}_m) = M^r_w$ with $r=\min\{1,p,q\}$. 
Clearly, the set $\{(ak,bj), k,j \in \Z^d\}$ is well-spread in $\HH^d$.
Thus, the assertion follows from Theorem \ref{thm_extend}.
\end{Proof} 

Depending on $p,q$ and $m$ one may also formulate weaker conditions on $g$ 
according to the abstract Theorem \ref{thm_extend}
such that the previous theorem still holds. 

Of course, one can also apply Theorem \ref{thm_nonlinapprox} to non-linear 
approximation with Gabor frames. We leave this  
straightforward task to the interested reader.
We only remark that apart from a small note in \cite{HoganLakey} 
modulation spaces with respect to Lorentz-spaces $L(p,q)$
did not appear in the literature in explicit form before. 

\section*{Acknowledgement}

Major parts of this paper were developed during a stay 
at the Mathematical Institute of the University of Wroc{\l}aw, Poland. The author 
would like to thank its members for their warm hospitality. 
In particular, he expresses his gratitude to M. Paluszy\'nski 
and R. Szwarc for interesting discussions on the subject. 
The author was supported 
by the European Union's Human Potential Programme under 
contracts HPRN-CT-2001-00273 (HARP) and  
HPRN-CT-2002-00285 (HASSIP).




\begin{thebibliography}{9}
\bibitem{DST} S. Dahlke, G. Steidl, and G. Teschke, Coorbit spaces and Banach frames on homogeneous spaces with applications to the sphere, Adv. Comp. Math., 21, 147--180, 2004. 
\bibitem{DST1} S. Dahlke, G. Steidl, and G. Teschke, Weighted coorbit spaces and Banach frames on homogeneous spaces, J. Four. Anal. Appl., {10}, 507--539, 2004.
\bibitem{DL} R.~DeVore, G.G.~Lorentz, {\it Constructive Approximation}, 
Springer, New York, 1993.
\bibitem{DM} M. Duflo, C.C. Moore, {On the regular representation of a nonunimodular locally compact group},
J. Funct. Anal., {21}, 209--243, 1976.
\bibitem{Fei_Wiener} {H.G. Feichtinger}, {Banach convolution algebras of Wiener type}, 
In: {\it Functions, Series, Operators, Proc. Int. Conf., Budapest 1980, Vol. I,} Colloq. Math. Soc. Janos Bolyai 35, pages 509--524, 1983. 
\bibitem{Fei_Homog} {H.G. Feichtinger}, {A characterization of minimal homogeneous Banach spaces}, {\it Proc. Am. Math. Soc.}, 81, 55--61, 1981.   
\bibitem{Fei_mod} {H.G. Feichtinger}, {Modulation spaces on locally compact groups},
Technical report, Vienna, 1983.
\bibitem{FG1} {H. G. Feichtinger, K. Gr\"ochenig}, {A unified approach to atomic decompositions via integrable group representations}, In: {\it Proc. Conference on Functions, Spaces and Applications, Lund, 1986}, Springer Lect. Notes Math., {1302}, 52--73, 1988.
\bibitem{FG2}  {H. G. Feichtinger, K. Gr\"ochenig}, {Banach spaces related to integrable group representations and their atomic decompositions I}, J. Funct. Anal., {86}, 307--340, 1989.
\bibitem{FG3} {H. G. Feichtinger, K. Gr\"ochenig}, {Banach spaces related to integrable group representations and their atomic decompositions II}, Monatsh. Math., {108}, 129--148, 1989.
\bibitem{FR} M. Fornasier, H. Rauhut, Continuous frames, function spaces, and the discretization problem, to appear in J. Fourier Anal. Appl., 2005.
\bibitem{GS} {Y.V. Galperin, S. Samarah}, 
{Time-frequency analysis on modulation spaces $M^{p,q}_m$, $0<p,q\leq \infty$}, 
Appl. Comp. Harm. Anal., 16, 1--18, 2004.
\bibitem{gro} K. Gr{\"o}chenig, Describing functions: atomic decompositions versus frames, Monatsh. Math., 112, 1--41, 1991.
\bibitem{Gr} {K. Gr\"ochenig}, {\it Foundations of Time-Frequency Analysis}, Birkh\"auser Verlag, 2001.
\bibitem{Gr3} {K. Gr\"ochenig}, {Unconditional bases in translation and dilation invariant
function spaces on $\R^n$}, In Sendov, B. et. al. (eds.), Constructive Theory of Functions, 
Proc. Conf. Varna 1987, 174--183, 1988.
\bibitem{GroS} {K. Gr\"ochenig, S. Samarah}, Nonlinear approximation with local Fourier bases, Constr. Approx. 16, 317--331, 2000.
\bibitem{HoganLakey} J.A. Hogan, J.D. Lakey, Sharp embeddings for modulation spaces and the Poisson summation formula, In:  Proc. Sampling Theory and Its Applications, 
Loen, Norway, 52--57, 1999.
\bibitem{Hunt} R.A. Hunt, On $L(p,q)$ spaces, Enseignement Math. 12, 249--276, 1966.
\bibitem{Kalton} N. Kalton, {Quasi-Banach spaces}, In: W.B. Johnson, J. Lindenstrauss (eds.), {\it Handbook of the Geometry of Banach spaces, Vol. 2}, North Holland, 1101--1129, 2003. 
\bibitem{hr} H. Rauhut, {\it Time-Frequency and Wavelet Analysis of Functions with Symmetry Properties}, Ph.D. thesis (2004), Technische Universit\"at M\"unchen, 
published at Logos-Verlag, 2005.
\bibitem{hr2} H. Rauhut, Banach frames in coorbit spaces consisting of elements which are invariant under symmetry groups, Appl. Comput. Harm. Anal., {18}, 94--122, 2005.
\bibitem{Rau_Wiener} H. Rauhut, Wiener amalgam spaces with respect 
to quasi-Banach spaces,
Preprint, 2005.
\bibitem{Triebel} {H. Triebel}, {Modulation spaces on the Euclidean $n$-space}, Z. Anal. Anwendungen 5, 443--457, 1983.
\end{thebibliography}
\end{document}